\documentclass[10pt]{article}
\usepackage{amsmath,amssymb,amsfonts, euscript}
\newtheorem{theorem}{Theorem}[section]

\begin{document}
\title{Invariance of the Barycentric Subdivision of a Simplicial Complex}
\author{Rashid Zaare-Nahandi}

\date{}

\maketitle

\begin{abstract}
In this paper we show that a simplicial complex can be determined uniquely up to isomorphism by its barycentric subdivision or comparability graph. At the end, it is summarized several algebraic, combinatorial and topological invariants of simplicial complexes.
\end{abstract}

\maketitle

\section{Introduction and Preliminaries}

Stanley-Reisner rings
of simplicial complexes, which have had fantastic application
in combinatorics [7], possess a rigidity property in the sense
that they determine their underlying simplicial complexes uniquely up to isomorphism [4, 8].
Barycentric subdivision of a simplicial complex is another very important and applicable object [1,2,3,7] which we want to prove that possess the same rigidity property.

Following, there are some basic definitions and facts on simplicial complexes and related topics which we will need later. See [3], [5] and [7] for details and more.

Let $[n] = \{1,2,\ldots,n\}$. A (finite) simplicial complex $\Delta$ on $n$ vertices, is a system of subsets of
$[n]$ such that the following conditions hold:\\
a) $\{i\}\in\Delta$ for any $i\in [n]$,\\
b) if $E\in\Delta$ and $F\subseteq E$, then $F\in\Delta$.\\
An element of $\Delta$ is called a face and a maximal face with respect to inclusion is called a facet. The set of all facets is denoted by $\mathcal{F}(\Delta)$. In the set of non-faces of $\Delta$ (those subsets of $[n]$ who are not in $\Delta$), the set of minimal elements is denoted by $\mathcal{N}(\Delta)$. The dimension of a face $F\in\Delta$ is defined to be $|F|-1$ and dimension of $\Delta$ is maximum of dimensions of its faces.

Let $\Delta$ be a simplicial complex on $[n]$ and  of dimension $d-1$. For each $0\leq i \leq d-1$ the $i$th skeleton of $\Delta$ is the simplicial complex $\Delta^{(i)}$ on $[n]$ whose faces
are those faces $F$ of $\Delta$ with $|F|\leq i + 1$. In particular the 1-skeleton $\Delta^{(1)}$ of
$\Delta$ is the finite graph on $[n]$ whose edges are the 1-dimensional faces $\{i, j\}$ of $\Delta$. We say that a simplicial complex $\Delta$ is connected if the finite graph $\Delta^{(1)}$ is connected.

Let $\Delta$ be a simplicial complex on
$[n]$. Let $S=K[x_1,\ldots,x_n]$ be the polynomial ring in $n$
indeterminates and with coefficients in a field $K$. Let $I_{\Delta}$
be the ideal of $S$ generated by all square free monomials
$x_{i_1}\cdots x_{i_s}$, provided that
$\{i_1,\ldots,i_s\}\not\in\Delta$. It is clear that the minimal generating set of $I_{\Delta}$ is all square free monomials $x_{i_1}\cdots x_{i_s}$, such that
$\{{i_1},\ldots,{i_s}\}\in\mathcal{N}(\Delta)$. The quotient ring
$K[\Delta]=S/I_{\Delta}$ is called the Stanley-Reisner ring of the simplicial
complex $\Delta$.

The facet ideal of $\Delta$ is the ideal $I(\Delta)$ of $S$ which is generated by those
square free monomials $x_{i_1}\cdots x_{i_s}$, provided that
$\{{i_1},\ldots,{i_s}\}$ is a facet in $\Delta$. The quotient ring $K_\mathcal{F}[\Delta]=S/I(\Delta)$ is called facet ring of $\Delta$.

For a given simplicial complex $\Delta$ on $[n]$, define $\Delta^{\vee}$ by
$$\Delta^{\vee} = \{[n] \setminus F : F \not\in \Delta\}.$$
It is clear that $\Delta^{\vee}$ is a simplicial complex and
$(\Delta^{\vee})^{\vee} = \Delta$. The simplicial complex $\Delta^{\vee}$ is called the Alexander dual of $\Delta$. Note that
$$\mathcal{F}(\Delta^{\vee}) = \{[n] \setminus F : F \in \mathcal{N}(\Delta)\}.\eqno{(1)}$$
Define the complement simplicial complex  $\Delta^c$ of $\Delta$ to be the simplicial complex whose facets are complements of facets of $\Delta$.  One has $$I_{\Delta^{\vee}} = I(\Delta^c)\eqno{(2)}$$.

A partially ordered set (poset) is a nonempty set $P$ with an order $\leq$ such that the followings hold\\
for each $x, y$ and $z$ in $P$, \\
a) $x\leq x$,\\
b) $x\leq y$ and $y\leq x$ implies $x=y$, \\
c) $x\leq y$ and $y\leq z$ implies $x\leq z$. \\
A simplicial complex can be assumed as a poset by order of inclusion.

Let $\Delta$ be a simplicial complex on vertex set $[n]$. The barycentric subdivision of $\Delta$, denoted by $\Delta^{\flat}$ is a simplicial complex which its vertex set is all elements of $\Delta$ other than the empty set, and two vertices are in a face if and only if one of them is subset of the other. In other words, facets of $\Delta^{\flat}$ are maximal chains in $\Delta$ assumed as a poset.

It is easy to check that the minimal non-faces of $\Delta^{\flat}$ are subsets of $\Delta$ with exactly two non-comparable elements. Therefore, $\Delta^{\flat}$ is a flag complex or a clique complex and the ideal $I_{\Delta^{\flat}}$ is generated by square-free quadrics. It is known that dimensions (and depths, respectively) of a simplicial complex and its barycentric subdivision are equal [6].

The 1-skeleton of $\Delta^{\flat}$ is called comparability graph of $\Delta$ and is denoted by $G(\Delta)$. The complement $\overline{G(\Delta)}$ of $G(\Delta)$ is called non-comparability graph of $\Delta$. The ideal $I_{\Delta^{\flat}}$ can be assumed as edge ideal of the graph $\overline{G(\Delta)}$, and then the simplicial complex $\Delta^{\flat}$ is the complex of independent sets of this graph.

It is not true that any graph is comparability graph of some simplicial complex. For example there is no any complex whose comparability graph is a cycle of length 3, 4 or 5. A necessary condition for a graph to be comparability graph of some complex is to be transitively orientable. That is, there is an orientation on the graph such that if $(x,y)$ and $(y,z)$ be oriented edges, the there is oriented edge $(x,z)$.

A (convex) polytope is the convex hull of a finite point set in Euclidean $n$ dimensional space, for some $n$.
A proper face of a polytope is the intersection of the polytope with
a supporting hyperplane. The empty set and the polytope itself are
called improper faces. A polyhedral complex is the union of a finite set of
polytopes such that, intersection of any
two being a face of each.

It is known that geometric realizations of $\Delta$ and $\Delta^{\flat}$ are homeomorphic as topological spaces and therefore, they share topological properties as Cohen-Macauleyness (see [7, p. 101]).

\section{The Main Result}

It is natural to ask whether a given graph is comparability graph of some simplicial complex and how many non isomorphic simplicial complexes are there with the same comparability graph. In this paper, we will prove that there is only one simplicial complex with a given comparability graph (up to isomorphism). Here, an isomorphism of simplicial complexes $\Delta_1$ and $\Delta_2$ is a bijection between their vertex sets who preserves faces and facets. It is enough to check that image and inverse image of any facet is again a facet. Face lattice of a polyhedral complex is a generalization of notion of simplicial complexes (see  [1] for definitions). In the case of polyhedral complexes with at least two maximal faces, M. Bayer has proved the following.

\begin{theorem} (Bayer [1]) Let $P$ be the face lattice of a connected polyhedral complex
with at least two maximal faces, and let $P^*$ be its dual poset. If $Q$
is a poset with $Q^{\flat} = P^{\flat}$, then either $Q = P$ or $Q = P^*$.
\end{theorem}

A simplicial complex is face lattice of the polyhedral complex of its geometric realization and so, by the above theorem, for a given connected simplicial cimplex with at least two facets, there are at most one other simplicial complex with the same barycentric subdivision. Let $\Delta$ be a simplicial complex, barycentric subdivision $\Delta^{\flat}$ is clique complex of the comparability graph $G(\Delta)$ and $G(\Delta)$ is 1-skeleton of $\Delta^{\flat}$. Therefore, knowing one of them is enough to construct the other. Therefore, two barycentric subdivisions are isomorphic as simplicial complexes if and only if their 1-skeletons are isomorphic as graphs. If $\Delta_1$ and $\Delta_2$ are two simplicial complexes, then they are isomorphic if and only if there is a  rearrangement of their connected components such that the corresponding components are isomorphic separately. Therefore, in the isomorphism problem of simplicial complexes, it is enough to consider only connected complexes, which is equivalent to have connected barycentric subdivision.

In the proof of the above theorem, M. Bayer has shown that, for a given $P$ with the mentioned conditions, there are exactly two transitive orientations on the graph of 1-skeleton of $P^{\flat}$, each reverse of the other. We will show that in case of simplicial complexes there is no need to the conditions and, either only one of the orientations on the graph admits a simplicial complex or, two simplicial complexes corresponding to these two orientations are isomorphic.

\begin{theorem}
Let $\Delta_1$ and $\Delta_2$ be simplicial complexes. Let two graphs $G(\Delta_1)$ and $G(\Delta_2)$ are isomorphic. Then, $\Delta_1$ and $\Delta_2$ are isomorphic.
\end{theorem}

{\it Proof.} As mentioned above, it is enough to consider connected simplicial complexes. Therefore, let $\Delta_1$ and $\Delta_2$ be connected. Without loss of generality we may assume that $G(\Delta_1) = G(\Delta_2) = G$.
If $\Delta_1$ is of dimension 0, then it is a single point and the theorem holds. Suppose that dimension of $\Delta_1=d\geq 1$. If $\Delta_1$ has only one facet, then it is a simplex and in the graph $G$, there is a unique vertex corresponding the facet, which is connected to all other vertices. The only possibility for $\Delta_2$ with a non empty face comparable with all others is to be a simplex. Therefore, $\Delta_2$ is a simplex of dimension $d$ and any two simplexes of the same dimension are isomorphic.

Now, we consider the case that $\Delta_1$ is not a simplex and has at least two facets. In this case, by the Theorem 2.1 of Bayer, there are at most two simplicial complexes with comparability graph $G$. One of them is $\Delta_1$. Let $\overrightarrow{G}$ be the orientation of $G$ corresponding to $\Delta_1$ such that $(x,y)$ is a directed edge if $x\subset y$ in $\Delta_1$. Give a grade to each vertex $x$ of $\overrightarrow{G}$ equal to the dimension of $x$ in $\Delta_1$. Note that, there is no any edge connecting two vertices with the same grade. The graph $\overrightarrow{G}$ is $(d+1)$-partite with each part consisting of all vertices of the same grade. In fact, there is a part consisting of all vertices which all arrows  are going out of them. This is the set of vertices of the underlying simplicial complex. All vertices of $\overrightarrow{G}$ such that all their connecting arrows are coming in, are facets and then the underlying simplicial complex can be uniquely determined. In the rest of proof, we investigate how a simplicial complex $\Delta_2$ may exist with comparability graph $G$ with revers orientation of $\overrightarrow{G}$, indicated by $\overleftarrow{G}$.

Suppose that there is such a simplicial complex $\Delta_2$. In a directed graph we call a vertex $v$ to be initial if there is no any edge connected to $v$ by direction with end point $v$. We call a vertex $v$ to be terminal if there is no any edge connected to it by direction with starting point $v$. In the above situation, terminal vertices of $\overrightarrow{G}$ are exactly initial vertices of $\overleftarrow{G}$ and vice versa.

First we show that $\overrightarrow{G}$ and $\overleftarrow{G}$ are "pure" in the sense that all maximal chains of them have a fixed length. In contrary, suppose there are two facets $F_1$ and $F_2$ in $\Delta_1$ with dim$(F_1)<$dim$(F_2)$. Denote their corresponding vertices in the graph $G$ also with the same names $F_1$ and $F_2$. These vertices are terminal points in $\overrightarrow{G}$ and initial points in $\overleftarrow{G}$. The complex $\Delta_2$ is connected, so, its 1-skeleton is connected. $F_1$ and $F_2$ are of dimension 0 in $\Delta_2$ and so there is no edge connecting them directly. Therefore, there is a vertex $E$ in $G$ with dimension one less than dimension of $F_1$ in $\Delta_1$ such that $F_1$ is connected directly to $E$ and it is connected directly to $F_2$ or another vertex with dimension equals to dim $F_2$. Therefore, $E$ is a proper subset of $F_2$ in $\Delta_1$
and so, it has an edge connecting to one of maximal proper subsets of $F_2$ as $E'$ which has dimension equal to  dim$F_2 - 1$. But, it is impossible because $E$ and $E'$ both have dimension 1 in $\Delta_2$ and can not have a common edge.

Now, suppose $\Delta_1$ is a pure simplicial complex with facets $F_1, F_2, \ldots, F_r$ and sub-facets (maximal proper subsets of facets) $E_1, E_2, \ldots, E_t$. The 1-skeleton of $\Delta_2$ consists of $F_1, F_2, \ldots, F_r$ as zero dimensional faces and $E_1$, $E_2$, $\ldots$, $E_t$ as one dimensional faces. Therefore, the induced subgraph of $G$
 to 1-skeleton of $\Delta_2$, consists of two parts  $F_1, F_2, \ldots, F_r$, all of the same degree and $E_1, E_2, \ldots, E_t$, all of degree 2. The induced graph has to be connected. Let $A=E_1\cup\cdots\cup E_t = F_1\cup\cdots\cup F_r$ where, $F_i$s and $E_j$s are considered as sets in $\Delta_1$. We claim that $F_1, F_2, \ldots, F_r$ are all maximal proper subsets of $A$. The claim is clearly true because only in this case we get a connected 2 partite graph with one part consisting of vertices of degree 2 and vertices of the other part all have the same degree. Therefore, the simplicial complex $\Delta_1$ is $(d-1)$-skeleton of a simplex. So is $\Delta_2$, and  because they have the same dimensions, they are isomorphic. \hfill $\Box$

Inverse of the above theorem is clearly true. That is, if $\Delta_1$ and $\Delta_2$ are isomorphic as simplicial complexes, then their barycentric subdivisions and comparability graphs will be isomorphic. Therefore, barycentric subdivision and comparability graph are invariants of a simplicial complex, and vise versa, underlying simpicial complex is invariant of its barycentric subdivision and comparability graph. In the following theorem we summarize more invariants of simplicial complexes.

\begin{theorem}
Let $\Delta_1$ and $\Delta_2$ be two simplicial complexes. The following conditions are equivalent.

\begin{enumerate}
  \item $\Delta_1$ and $\Delta_2$ are isomorphic as simplicial complexes.
  \item $\Delta_1^{\vee}$ and $\Delta_2^{\vee}$ are isomorphic as simplicial complexes.
  \item $\Delta_1^c$ and $\Delta_2^c$ are isomorphic as simplicial complexes.
  \item $\Delta_1^{\flat}$ and $\Delta_2^{\flat}$ are isomorphic as simplicial complexes.
  \item $\Delta_1^{\flat^n}$ and $\Delta_2^{\flat^n}$ are isomorphic as simplicial complexes for some positive integer $n$.
  \item $K[\Delta_1]$ and $K[\Delta_2]$ are isomorphic as $K$ algebras.
  \item $K_{\mathcal{F}}[\Delta_1]$ and $K_{\mathcal{F}}[\Delta_2]$ are isomorphic as $K$ algebras.
  \item $G(\Delta_1)$ and $G(\Delta_2)$ are isomorphic as graphs.
  \end{enumerate}
\end{theorem}

{\it Proof.} Equivalences of 1, 2 and 3 are clear form equations (1) and (2). Equivalence of 1 and 4 is proved in Theorem 2.2 above. Items 4 and 5 are equal because of 1 and 4. Equivalences of 1, 6 and 7 are proved in [4] and [8]. By the  argument just after Theorem 2.1, equivalence of 4 and 8 is clear. \hfill $\Box$

\noindent
Rashid Zaare-Nahandi\\
Institute for Advanced Studies in Basic Sciences\\
P.O.Box 45195-1159\\
Zanjan, Iran\\
rashidzn@iasbs.ac.ir

\end{document}